\documentclass[12pt]{article}%
\usepackage{amsmath}
\usepackage{amsfonts}
\usepackage{amssymb}
\usepackage{ytableau}
\usepackage{graphicx}
\providecommand{\U}[1]{\protect\rule{.1in}{.1in}}
\newtheorem{theorem}{Theorem}

\begin{document}

\title{5-Engel Lie algebras}
\author{Michael Vaughan-Lee}
\date{January 2024}
\maketitle

\begin{abstract}
We describe computer calculations which show that if $L$ is a 5-Engel Lie
algebra over a field of characteristic zero, or over a field of prime
characteristic $p>7$,$\,$\ then $L$ is nilpotent of class at most 11. We use
the representation theory of the symmetric group to show that the problem can
be reduced to showing that certain four generator Lie superalgebras satisfying
relations derived from the 5-Engel identity are nilpotent of class at most 11.
We also describe computer calculations which show that if $G$ is a finite
5-Engel $p$-group for a prime $p>7$ then $G$ is nilpotent of class at most 10.

\end{abstract}

\section{Introduction}

A Lie algebra $L$ is said to be an $n$-Engel Lie algebra if ad$(x)^{n}=0$ for
all $x\in L$. Clearly 1-Engel Lie algebras are abelian, and Levi \cite{levi}
showed that 2-Engel Lie algebras are nilpotent of class at most 3. Heineken
\cite{heineken61} proved that 3-Engel Lie algebras over a field of
characteristic zero, or over a field of prime characteristic $p\neq2,5$, are
nilpotent of class at most 4. And Traustason \cite{Traustason93} proved that
4-Engel Lie algebras over a field of characteristic zero, or over a field of
prime characteristic $p>5$, are nilpotent of class at most 7. In this note we
prove the following theorem.

\begin{theorem}
If $L$ is a $5$-Engel Lie algebra over a field of characteristic zero, or over
a field of prime characteristic $p>7$, then $L$ is nilpotent of class at most
$11$.
\end{theorem}

It follows from Zel'manov's solution of the restricted Burnside problem
\cite{zelmanov91}, \cite{zelmanov91b} that $n$-Engel Lie algebras are locally
nilpotent (for all $n$). And in a very deep theorem Zel'manov
\cite{zelmanov87} also showed that $n$-Engel Lie algebras over a field of
characteristic zero are (globally) nilpotent (for all $n$). Note that if
$n$-Engel Lie algebras over fields of characteristic zero are nilpotent of
class at most $c$ then $n$-Engel Lie algebras of prime characteristic $p$ will
also be nilpotent of class at most $c$ for all sufficiently large $p$.\ On the
other hand Razmyslov \cite{Razmyslov71} has shown that there are non-soluble
(and hence non-nilpotent) $(p-2)$-Engel Lie algebras of characteristic $p$ for
all primes $p\geq5$. For $p=5$ this result was originally proved by Bachmuth
and Mochizuki \cite{bachmuth70}.

It is easy to construct non-nilpotent 3-Engel Lie algebras of characteristic
2. Let $L$ be the free metabelian Lie algebra over $\mathbb{Z}_{2}$, freely
generated by $x_{1},x_{2},\ldots$. Then $[L,L]$ is an abelian Lie algebra with
a vector space basis consisting of all left normed Lie products%
\[
\lbrack x_{i_{1}},x_{i_{2}},x_{i_{3}},\ldots,x_{i_{n}}]
\]
with $n\geq2$, $i_{1}>i_{2}\leq i_{3}\leq\ldots\leq i_{n}$. We let $I$ be the
ideal of $L$ generated by all such Lie products where $i_{j}=i_{k}$ with
$1\leq j<k\leq n$. Then $L/I$ is non-nilpotent since%
\[
\lbrack x_{2},x_{1},x_{3},x_{4},\ldots,x_{n}]\notin I.
\]
On the other hand, it is easy to see that $L/I$ satisfies the identical
relation $[x,y,z,z]=0$, so that $L/I$ satisfies the 3-Engel identity.

The same construction over the field $\mathbb{Z}_{3}$ gives a non-nilpotent
Lie algebra satisfying the 4-Engel identity.

A group is said to be an $n$-Engel group if it satisfies the identity
$[x,\underset{n}{\underbrace{y,y,\ldots,y}}]=1$. Lie algebra methods can be
applied to finite $n$-Engel $p$-groups, and Levi \cite{levi} shows that
2-Engel groups are nilpotent of class at most 3. Similarly Heineken
\cite{heineken61} proves that 3-Engel groups are locally nilpotent, and are
nilpotent of class at most 4 provided they have no elements of order 2 or 5.
Traustason \cite{Traustason93} proves that locally nilpotent 4-Engel
$p$-groups are nilpotent of class at most 7 provided $p>5$. (Note that by a
result of Havas and Vaughan-Lee \cite{havasvl05}\ we know that all 4-Engel
groups are locally nilpotent.) In this note we prove the following theorem.

\begin{theorem}
If $G$ is a locally nilpotent $5$-Engel $p$-group for some prime $p>7$ then
$G$ is nilpotent of class at most $10$.
\end{theorem}

\section{Free Lie algebras}

Let $L$ be a free Lie algebra over a commutative ring $R$ with 1, and let $L$
have free generating set $A$. Then $L$ is a free module over $R$ and has a
basis as a free module consisting of basic commutators (or basic Lie products).

The basic commutators of weight one are the elements $a\in A$, which we assume
to be an ordered set.

The basic commutators of weight two are the elements $[a,b]$ where $a,b\in A$
and $a>b$. These are ordered arbitrarily among themselves, and so that they
follow basic commutators of weight one.

The basic commutators of weight three are the elements $[a,b,c]$ where
$a,b,c\in A$ and $a>b\leq c$. (We use the left-normed convention so that
$[a,b,c]$ denotes $[[a,b],c]$.) The basic commutators of weight three are
ordered arbitrarily among themselves, and so that they follow the basic
commutators of weight two.

In general, if $k>3$ the basic commutators of weight $k$ are the commutators
$[c,d]$ where for some $m,n$ such that $m+n=k$,

\begin{enumerate}
\item $c,d$ are basic commutators of weight $m,n$ respectively,

\item $c>d$,

\item if, in the definition of basic commutators of weight $m$, $c$ was
defined to be $[e,f]$ then $f\leq d$.
\end{enumerate}

The basic commutators of weight $k$ are then ordered arbitrarily among
themselves, and so that they follow the basic commutators of weight $k-1$.

A proof that the basic commutators form an $R$-module basis for $L$ can be
found in my book \cite{vlee93b}.

\section{Lie superalgebras}

A Lie superalgebra is a $\mathbb{Z}_{2}$-graded algebra $L=L_{0}\oplus L_{1}$
with a bilinear product $[,]$ such that
\begin{align*}
\lbrack L_{0},L_{0}],\,[L_{1},L_{1}] &  \leq L_{0},\\
\lbrack L_{0},L_{1}],\,[L_{1},L_{0}] &  \leq L_{1}.
\end{align*}
Elements in $L_{0}$ are said to be even elements, and elements in $L_{1}$ are
said to be odd. If $a$ is even then we set $|a|=0$, and if $a$ is odd then we
set $|a|=1$. Odd elements and even elements are said to be homogeneous.
Finally, the product $[,]$ must satisfy the following relations for all
homogeneous elements $a,b,c$.%
\[
\lbrack b,a]=-(-1)^{|a|.|b|}[a,b].
\]%
\[
(-1)^{|a|.|c|}[a,[b,c]]+(-1)^{|b|.|a|}[b,[c,a]]+(-1)^{|c|.|b|}[c,[a,b]]=0.
\]
It is helpful to note that these relations imply that if $a,b,c$ are
homogeneous elements then%
\[
\lbrack a,[b,c]]=[a,b,c]-(-1)^{|b|.|c|}[a,c,b].
\]
We also add in the requirement that $[a,a]=0$ for even elements, and the
requirement that $[a,a,a]=0$ for odd elements. (These extra requirements are
redundant if 2 and 3 are invertible.)

Now let $L$ be a free Lie superalgebra over a commutative ring $R$ with 1, and
let $L$ have free basis $A$, where some of the elements of $A$ are even, and
the rest are odd. Let $C$ be a set of basic commutators on $A$. (Note that the
definition of basic commutators is not affected by the fact that some of the
elements of $A$ are even and some are odd. Also note that the elements of $C$
are homogeneous.) There is a proof of the following theorem in \cite{arxiv}.

\begin{theorem}
The free Lie superalgebra $L$ is a free $R$-module with basis consisting of $C$
together with Lie products $[c,c]$ where $c \in C$ is odd.
\end{theorem}

We will make use of the following construction for Lie superalgebras. We let
$G$ be the associative algebra over $R$ generated by $e_{1},e_{2},\ldots$
subject to the relations $e_{i}^{2}=0$ for $i=1,2,\ldots$, and $e_{i}%
e_{j}=-e_{j}e_{i}$ for all $i,j$ with $i\neq j$. So $G$ has an $R$-module
basis $B$ consisting of all possible products $e_{i}e_{j}\ldots e_{k} $ with
$i<j<\ldots<k$. We can write $G=G_{0}\oplus G_{1}$ where $G_{0} $ is spanned
by the products $e_{i}e_{j}\ldots e_{k}$ of even length and $G_{1}$ is spanned
by the products $e_{i}e_{j}\ldots e_{k}$ of odd length. If $g\in G_{0}$ and
$h\in G$ then $gh=hg$, but if $g,h\in G_{1}$ then $gh=-hg$. Now let $K$ be a
Lie algebra over $R$ and let $L=K\otimes G$. We define a bracket product on
$L$ by setting%
\[
\lbrack a\otimes g,b\otimes h]=[a,b]\otimes gh
\]
for $a,b\in K$ and $g,h\in G$, and extending this product to the whole of $L$
by linearity. If we set $L_{0}=K\otimes G_{0}$, $L_{1}=K\otimes G_{1}$ then
$L=L_{0}\oplus L_{1}$ and it is straightforward to show that $L$ is a Lie superalgebra.

We now suppose that $K$ satisfies the 5-Engel identity, and deduce some
relations which hold in $L$. Let $a,a_{1},a_{2},a_{3},a_{4},a_{5}\in K$. Then%
\begin{equation}
\sum_{\sigma\in\text{Sym}(5)}[a,a_{1\sigma},a_{2\sigma},a_{3\sigma}%
,a_{4\sigma},a_{5\sigma}]=0.
\end{equation}
Now let $g,g_{1},g_{2},g_{3},g_{4},g_{5}$ be elements of the basis $B$ for $G
$, and consider%
\[
\lbrack a\otimes g,a_{1\sigma}\otimes g_{1\sigma},a_{2\sigma}\otimes
g_{2\sigma},a_{3\sigma}\otimes g_{3\sigma},a_{4\sigma}\otimes g_{4\sigma
},a_{5\sigma}\otimes g_{5\sigma}]
\]
for $\sigma\in $Sym(5). Expanding we obtain%
\[
\lbrack a,a_{1\sigma},a_{2\sigma},a_{3\sigma},a_{4\sigma},a_{5\sigma}]\otimes
gg_{1\sigma}g_{2\sigma}g_{3\sigma}g_{4\sigma}g_{5\sigma}=\varepsilon\lbrack
a,a_{1\sigma},a_{2\sigma},a_{3\sigma},a_{4\sigma},a_{5\sigma}]\otimes
gg_{1}g_{2}g_{3}g_{4}g_{5}%
\]
where $\varepsilon=\pm1$. We can compute $\varepsilon$ as follows. We take the
sequence $(g_{1},g_{2},g_{3},g_{4},g_{5})$ and delete all the even elements
from this sequence leaving a sequence $(g_{i},g_{j},\ldots,g_{k})$ where
$g_{i},g_{j},\ldots,g_{k}$ are odd. Similarly we delete the even elements from
the sequence $(g_{1\sigma},g_{2\sigma},g_{3\sigma},g_{4\sigma},g_{5\sigma})$
leaving the sequence $(g_r,g_s,\ldots,g_t)$ where $(r,s,\ldots ,t)$ is a
permutation of $(i,j,\ldots ,k)$. Then
$\varepsilon=1$ if $(r,s,\ldots ,t)$ is an even
permutation of $(i,j,\ldots ,k)$ and otherwise $\varepsilon=-1$. So
$\varepsilon$ depends on the permutation $\sigma$ induces on the odd elements
in the sequence $(g_{1},g_{2},g_{3},g_{4},g_{5})$. We use the notation
$|\sigma_{\text{odd}}| $ to denote the value of $\varepsilon$. So%
\begin{align*}
& [a\otimes g,a_{1\sigma}\otimes g_{1\sigma},a_{2\sigma}\otimes g_{2\sigma
},a_{3\sigma}\otimes g_{3\sigma},a_{4\sigma}\otimes g_{4\sigma},a_{5\sigma
}\otimes g_{5\sigma}]\\
& =|\sigma_{\text{odd}}|[a,a_{1\sigma},a_{2\sigma},a_{3\sigma},a_{4\sigma
},a_{5\sigma}]\otimes gg_{1}g_{2}g_{3}g_{4}g_{5}.
\end{align*}
(Note that $gg_{1}g_{2}g_{3}g_{4}g_{5}=0$ if and only if
$gg_{1\sigma}g_{2\sigma}g_{3\sigma}g_{4\sigma}g_{5\sigma}=0$
for all permutations $\sigma $.)
It follows from equation (1) that
\[
\sum_{\sigma\in\text{Sym(5)}}|\sigma_{\text{odd}}|[a\otimes g,a_{1\sigma
}\otimes g_{1\sigma},a_{2\sigma}\otimes g_{2\sigma},a_{3\sigma}\otimes
g_{3\sigma},a_{4\sigma}\otimes g_{4\sigma},a_{5\sigma}\otimes g_{5\sigma}]=0.
\]
It follows by linearity that if $b,b_{1},b_{2},b_{3},b_{4},b_{5}$ are
homogeneous elements of $L$ then%
\[
\sum_{\sigma\in\text{Sym(5)}}|\sigma_{\text{odd}}|[b,b_{1\sigma},b_{2\sigma
},b_{3\sigma},b_{4\sigma},b_{5\sigma}]=0,
\]
where $|\sigma_{\text{odd}}|$ is the sign of the permutation $\sigma$ induces
on the odd elements of the sequence $(b_{1},b_{2},b_{3},b_{4},b_{5})$.

\section{Representation theory of the symmetric group}

We let $N$ be a positive integer, and we consider the group ring $\mathbb{Q}%
$Sym$(N)$ of the symmetric group on $N$ letters, where $\mathbb{Q}$ is the
rational field. The identity element in $\mathbb{Q}$Sym$(N)$ is a sum of
primitive idempotents, and these are described in James and Kerber
\cite{jamesker81}: they correspond to \emph{Young tableaux}. For each
partition $(m_{1},m_{2},\ldots,m_{s})$ of $N$ with $m_{1}\geq m_{2}\geq
\ldots\geq m_{s} $ we associate a \emph{Young diagram}, which is an array of
$N$ boxes arranged in $s$ rows, with $m_{i}$ boxes in the $i$-th row. The
boxes are arranged so that the $j$-th column of the array consists of the
$j$-th boxes out of the rows which have length $j$ or more. For example, if
$N=5$ there are seven possible Young diagrams.

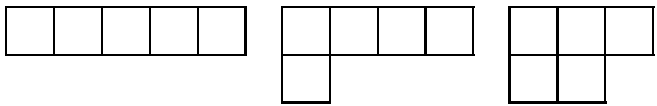
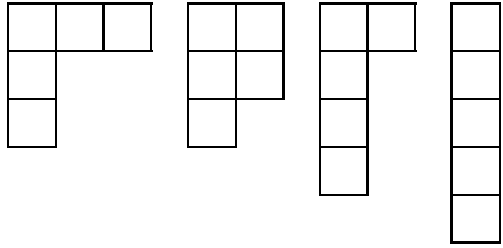
\[
\ydiagram{5}\;\;\;\;
\ydiagram{4,1}\;\;\;\;
\ydiagram{3,2}
\]
\bigskip

\[
\ydiagram{3,1,1}\;\;\;\;
\ydiagram{2,2,1}\;\;\;\;
\ydiagram{2,1,1,1}\;\;\;\;
\ydiagram{1,1,1,1,1}
\]
\bigskip

We obtain a Young tableau from a Young diagram by filling in the $N$ boxes
with $1,2,\ldots,N$ in some order. We then let $H$ be the subgroup of Sym$(N)$
which permutes the entries within each row of the tableau, and we let $V $ be
the subgroup of Sym$(N)$ which permutes the entries within each column of the
tableau. We set%
\[
e=\sum_{\pi\in V,\,\rho\in H}\text{sign}(\pi)\pi\rho.
\]
Then $\frac{1}{k}e$ is a primitive idempotent of $\mathbb{Q}$Sym$(N)$ for some
$k$ dividing $N!$. As mentioned above the identity element in $\mathbb{Q}%
$Sym$(N)$ can be written as a sum of primitive idempotents of this form, and
if the field $F$ has characteristic zero, or characteristic coprime to $N!$,
then so can the identity element in $F$Sym$(N)$.

Now let $K$ be the free 5-Engel Lie algebra generated by $a_{1},a_{2}%
,\ldots,a_{12}$ over a field $F$ where $F$ has characteristic zero, or prime
characteristic $p>12$. We want to prove that $[a_{1},a_{2},\ldots,a_{12}]=0$.
Let $M$ be the subspace of $K$ spanned by all commutators $[a_{1\sigma
},a_{2\sigma},\ldots,a_{12\sigma}]$ where $\sigma\in\,$Sym(12). There is a
natural action of $F$Sym(12) on $M$ obtained by setting%
\[
\lbrack a_{1\sigma},a_{2\sigma},\ldots,a_{12\sigma}]\tau=[a_{1\sigma\tau
},a_{2\sigma\tau},\ldots,a_{12\sigma\tau}]
\]
for $\tau\in\,$Sym(12). To prove that $[a_{1},a_{2},\ldots,a_{12}]=0$ it is
sufficient to show that
\[
[a_{1},a_{2},\ldots,a_{12}]e=0
\]
for all primitive idempotents in $F$Sym(12), as described above.

Consider a Young diagram with 12 cells. It is easy to see that it can be
subdivided into at most 4 disjoint horizontal and vertical strips. If the
first column of the diagram has length $k$ then we can divide the Young
diagram into $k$ horizontal strips. So we need only consider the case when
$k>4$. If we remove the first column we are left with a Young diagram on $12-k
$ cells. If the second column of the diagram has length $m$ then we can divide
up the diagram into one vertical strip of length $k$ and $m$ horizontal
strips. So we may assume that $m>3$. So there are at most 3 cells outside the
first two columns, and these cells must lie in a single vertical strip, or in
a single horizontal strip, or in two horizontal strips of lengths 2 and 1.

Suppose for example that we have a Young tableau with rows of length 5,3,2,1,1.
Then it can be divided up into one vertical strip of length 5, and three
horizontal strips of lengths 4, 2, 1. Suppose that the entries in the vertical
strip are $i_{1},i_{2},i_{3},i_{4},i_{5}$ and that the entries in the three
horizontal strips are $i_{6},i_{7},i_{8},i_{9}$ and $i_{10},i_{11}$ and
$i_{12}$. Let $T_{1}$ be the symmetric group on $\{i_{1},i_{2},i_{3}%
,i_{4},i_{5}\}$, and let $T_{2},T_{3},T_{4}$ be the symmetric groups on
$\{i_{6},i_{7},i_{8},i_{9}\}$, $\{i_{10},i_{11}\}$, $\{i_{12}\}$. Let $V$ be
the subgroup of Sym(12) which permutes the entries within each column of the
tableau, and let $H$ be the subgroup of Sym(12) which permutes the entries
within each row of the tableau. So $T_{1}$ is a subgroup of $V$, and
$T_{2}\times T_{3}\times T_{4}$ is a subgroup of $H$, and we can write%
\[
e=\sum_{\pi\in V,\,\rho\in H}\text{sign}(\pi)\pi\rho.
\]
as a sum of terms%
\[
\sum_{\sigma,\tau}\text{sign}(\sigma)\sigma\left(  \sum_{\pi\in T_{1},\,\rho\in
T_{2}\times T_{3}\times T_{4}}\text{sign}(\pi)\pi\rho\right)\tau
\]
where $\sigma$ runs over a left transversal for $T_{1}$ in $V$ and where
$\tau$ runs over a right transversal for $T_{2}\times T_{3}\times T_{4}$ in
$H$. So to show that%
\[
\lbrack a_{1},a_{2},\ldots,a_{12}]e=0
\]
it is sufficient to show that
\[
\lbrack a_{1},a_{2},\ldots,a_{12}]\sigma\left(  \sum_{\pi\in T_{1},\,\rho\in
T_{2}\times T_{3}\times T_{4}}\text{sign}(\pi)\pi\rho\right)  =0
\]
for all $\sigma$ in a left transversal for $T_{1}$ in $V$. We let $L=K\otimes
G$ as described above, and let%
\begin{align*}
b_{1}  & =a_{i_{1}}\otimes e_{1}+a_{i_{2}}\otimes e_{2}+a_{i_{3}}\otimes
e_{3}+a_{i_{4}}\otimes e_{4}+a_{i_{5}}\otimes e_{5},\\
b_{2}  & =a_{i_{6}}\otimes e_{6}e_{7}+a_{i_{7}}\otimes e_{8}e_{9}+a_{i_{8}%
}\otimes e_{10}e_{11}+a_{i_{9}}\otimes e_{12}e_{13},\\
b_{3}  & =a_{i_{10}}\otimes e_{14}e_{15}+a_{i_{11}}\otimes e_{16}e_{17},\\
b_{4}  & =a_{i_{12}}\otimes e_{18}e_{19}.
\end{align*}
Then%
\begin{align*}
& [a_{1},a_{2},\ldots,a_{12}]\sigma\left(  \sum_{\pi\in T_{1},\,\rho\in
T_{2}\times T_{3}\times T_{4}}\text{sign}(\pi)\pi\rho\right)  \otimes
e_{1}e_{2}\ldots e_{19}\\
& =\pm\lbrack c_{1},c_{2},\ldots,c_{12}]
\end{align*}
where five of the entries $c_{i}$ are equal to $b_{1}$, four of the entries
are equal to $b_{2}$, two of the entries are equal to $b_{3}$ and one of the
entries is equal to $b_{4}$. To be precise,%
\[
\lbrack a_{1},a_{2},\ldots,a_{12}]\sigma=[a_{1\sigma},a_{2\sigma}%
,\ldots,a_{12\sigma}],
\]
and if $\{i\sigma,j\sigma,k\sigma,l\sigma,m\sigma\}=\{i_{1},i_{2},i_{3}%
,i_{4},i_{5}\}$, then $c_{i}=c_{j}=c_{k}=c_{l}=c_{m}=b_{1}$, and so on. So to
show that $[a_{1},a_{2},\ldots,a_{12}]e=0$ it is sufficient to show that any
Lie product of weight 12 in $L$ is zero if it has weight 5 in $b_{1}$, weight
4 in $b_{2}$, weight 2 in $b_{3}$ and weight 1 in $b_{4}$. Clearly to prove
this we only need to consider the subalgebra of $L$ generated by $b_{1}%
,b_{2},b_{3},b_{4}$. So it is sufficient to show that if $L$ is a Lie
superalgebra generated by one odd element and three even elements, and if $L$
satisfies relations derived from the 5-Engel identity as described in Section
3, then Lie products of multiweight (5,4,2,1) in the four generators are zero.
(We will describe how to accomplish this in Section 5 below.)

Next consider a Young tableau with rows of length 4,3,2,2,1. The first column
has length five and the second column has length four. Removing these two
columns we are left with two horizontal strips of length two and length one.
Let the entries in the cells of the first two columns be $i_{1},i_{2}%
,i_{3},i_{4},i_{5}$ and $i_{6},i_{7},i_{8},i_{9}$, and let the entries in the
two remaining horizontal strips be $i_{10},i_{11}$ and $i_{12} $. As above we
let $T_{1},T_{2},T_{3},T_{4}$ be the symmetric groups on the sets
$\{i_{1},i_{2},i_{3},i_{4},i_{5}\}$, $\{i_{6},i_{7},i_{8},i_{9}\} $,
$\{i_{10},i_{11}\}$, $\{i_{12}\}$. Let $V$ be the subgroup of Sym(12) which
permutes the entries within each column of the tableau, and let $H$ be the
subgroup of Sym(12) which permutes the entries within each row of the tableau.
So $T_{1}\times T_{2}$ is a subgroup of $V$, and $T_{3}\times T_{4}$ is a
subgroup of $H$, and we can write%
\[
e=\sum_{\pi\in V,\,\rho\in H}\text{sign}(\pi)\pi\rho.
\]
as a sum of terms%
\[
\sum_{\sigma,\tau}\text{sign}(\sigma)\sigma\left(  \sum_{\pi\in T_{1}\times
T_{2},\,\rho\in T_{3}\times T_{4}}\text{sign}(\pi)\pi\rho\right)\tau
\]
where $\sigma$ runs over a left transversal for $T_{1}\times T_{2}$ in $V$ and
where $\tau$ runs over a right transversal for $T_{3}\times T_{4}$ in $H $. So
to show that%
\[
\lbrack a_{1},a_{2},\ldots,a_{12}]e=0
\]
it is sufficient to show that
\[
\lbrack a_{1},a_{2},\ldots,a_{12}]\sigma\left(  \sum_{\pi\in T_{1}\times T_{2},\,\rho\in
 T_{3}\times T_{4}}\text{sign}(\pi)\pi\rho\right)  =0
\]
for all $\sigma$ in a left transversal for $T_{1}\times T_{2}$ in $V$. We let $L=K\otimes
G$, and let%
\begin{align*}
b_{1}  & =a_{i_{1}}\otimes e_{1}+a_{i_{2}}\otimes e_{2}+a_{i_{3}}\otimes
e_{3}+a_{i_{4}}\otimes e_{4}+a_{i_{5}}\otimes e_{5},\\
b_{2}  & =a_{i_{6}}\otimes e_{6}+a_{i_{7}}\otimes e_{7}+a_{i_{8}}\otimes
e_{8}+a_{i_{9}}\otimes e_{9},\\
b_{3}  & =a_{i_{10}}\otimes e_{10}e_{11}+a_{i_{11}}\otimes e_{12}e_{13},\\
b_{4}  & =a_{i_{12}}\otimes e_{14}e_{15}.
\end{align*}
Then for each $\sigma$%
\[
\lbrack a_{1},a_{2},\ldots,a_{12}]\sigma\left(  \sum_{\pi\in T_{1}\times
T_{2},\,\rho\in T_{3}\times T_{4}}\text{sign}(\pi)\pi\rho\right)  \otimes
e_{1}e_{2}\ldots e_{15}%
\]
is a product of multiweight (5,4,2,1) in the generators $b_{1},b_{2}%
,b_{3},b_{4}$. So to show that $[a_{1},a_{2},\ldots,a_{12}]e=0$ it is
sufficient to show that if $L$ is a Lie superalgebra generated by two odd
element and two even elements, and if $L$ satisfies relations derived from the
5-Engel identity, then Lie products of multiweight (5,4,2,1) in the four
generators are zero.

In this way we can reduce the problem of showing that $[a_{1},a_{2}%
,\ldots,a_{12}]=0$ to a number of calculations in three and four generator Lie
superalgebras. (The actual number of calculations turns out to be 42, together
with some subsidiary calculations needed to handle three particularly difficult cases.)

\subsection{Characteristic 11}

We cannot use the representation theory of Sym(12) in characteristic 11, but
we can use the representation of Sym(8) to show that in a 5-Engel Lie algebra
of characteristic 11 all products of multiweight (4,1,1,1,1,1,1,1,1) in 9
generators are zero. As above, the representation theory of Sym(8) can be used
to reduce this problem to a number of calculations in three and four generator
Lie superalgebras. (Many of the calculations are covered by the calculations
for general characteristic $p>11$.) We can then use the representation theory
of Sym(10) to show that in a 5-Engel Lie algebra of characteristic 11 all
products of multiweight (2,1,1,1,1,1,1,1,1,1,1) in 11 generators are zero.
Again we can use the representation theory to reduce the problem to a number
of calculations in three and four generator superalgebras. Since we can assume
that all products of multiweight (4,1,1,1,1,1,1,1,1) are zero we only need
consider Young tableaux on 10 cells where the first row has length at most 3.

All this shows that 5-Engel Lie algebras of characteristic 11 satisfy the
identical relation
\[
[y,x,x,x_1,x_2,x_3,x_4,x_5,x_6,x_7,x_8,x_9]=0
\]
and this implies that they are nilpotent of class at most 11.

\section{The calculations}

We need to calculate in three and four generator Lie superalgebras satisfying
relations derived from the 5-Engel identity. For each of these
algebras we need to show that products with certain multiweights in the
generators are zero. It is important to note that when computing a three
or four generator Lie superalgebra we only need relations derived
from the 5-Engel identity that are homogenous in every generator. This means that
if (for example) we are calculating products of multiweight (5,4,2,1) in
generators $b_{1},b_{2},b_{3},b_{4}$, then we can assume throughout the
calculation that any Lie product with 6 or more entries $b_{1}$ is zero, that any
product with 5 or more entries $b_{2}$ is zero, and so on. Without this reduction most
of the calculations would be infeasible.

For characteristic $p<2^{16}$ these calculations can all be carried out using the
nilpotent quotient algorithm for graded Lie rings (see \cite{havasnvl90}).
With minor modifications the program is able to compute nilpotent quotients of
finitely generated Lie superalgebras over GF$(p)$ satisfying a finite set of
homogeneous relations. (The program can also be modified to handle larger
primes, with the limitation that the primes have to be stored in a 32 bit
word.) The program first computes the class one quotient of the algebra, then
computes the class two cover of this algebra, and enforces any weight two
relations. This gives the class two quotient of the algebra. Then the program
computes the class 3 cover of the class two quotient, and enforces the weight
3 relations. And so on.

Carrying out these calculations using the nilpotent quotient algorithm for
primes 11, 13, 17, 19 provides evidence that Theorem 1 holds true for all
primes $p>7$, and that it also holds true in characteristic zero. 
But to prove Theorem 1 for infinitely many primes we need a different approach. 

I wrote a nilpotent quotient program in
\textsc{Magma} \cite{boscan95} to compute Lie superalgebras over the rationals.
Unfortunately enforcing the 5-Engel relations class by class proved to be
infeasible, as there was a combinatorial explosion in the denominators of
the structure constants in the algebra. When computing free Lie superalgebras
over the rationals the structure constants were manageable, but they were
not all integers. The $p$-quotient algorithm in groups and the nilpotent quotient 
algorithm in Lie algebras described in \cite{havasnvl90} builds up a basis at each
class consisting of left-normed commutators (or Lie products), and
it \emph{may} be possible in a free Lie superalgebra to choose a basis
of left-normed Lie products which give integer structure constants.
But I was unable to find a way of doing this.

So I wrote a
program to compute free Lie superalgebras over the integers using the basis
described in Theorem 3 above consisting of basic commutators and Lie products
$[c,c]$ where $c$ is an odd basic commutator. When calculating products of
multiweight (5,4,2,1) in generators $b_{1},b_{2},b_{3},b_{4}$ (for example) I
was able to add in the relations that any basic commutator with 6 or more entries
$b_{1}$ was zero, and so on. So in every case I ended up with a basis for a
class 12 Lie superalgebra where all the structure constants giving the Lie
product of two basis elements were integers. I then computed all relations of
weight 12 derived from the 5-Engel identity and stored them in an integer
matrix. In most cases I was able to compute the elementary divisors of this
relation matrix, and check that the rank of the matrix was equal to the number
of columns, and that the elementary divisors only involved the
primes 2,3,5,7. This showed that in characteristic $p>7$, or in
characteristic zero, the relations implied that all products of weight 12 were zero.

This worked well for most of the relation matrices, but my computer with 16 gb
of RAM ran out of space for some of the bigger ones. (The smallest relation
matrix had size 277$\times$110, but one of the bigger ones had size
15791$\times$8400, and I was unable to compute the elementary divisors of this matrix.)

So I tried another approach. Suppose we have an $m\times n$ integer matrix
$(m>n)$, and suppose we want to prove that it has rank $n$ when read as a
matrix over GF$\left(  p\right)  $ for $p>7$. (This will also show that it has
rank $n$ when viewed as a matrix over the rationals.) Randomly pick $n$ rows
which are linearly independent when read as vectors over GF(11). Compute the
determinant of the $n\times n$ matrix with these rows, and divide out the
determinant by powers of 2, 3, 5, 7. Store the result. Repeat this with
another one or two random $n\times n$ matrices, and compute the greatest
common divisor of the two or three stored results. In most cases the greatest
common divisor of two of the results was 1, but in a very few cases I needed
to take the greatest common divisor of three of these values to get greatest
common divisor 1. This was enough to show that the matrix had rank $n$ when
thought of as a matrix over GF$(p)$ for $p>7$. This method used much less 
computer memory and was much faster than computing the 
elementary divisors of the matrix.

There was one particular calculation where this method gave a very
extraordinary result. Calculations with the nilpotent quotient algorithm show
that four generator 5-Engel Lie algebras (all four generators even) have class
10 rather than class 11 for $p=11,13,17,19$. It seemed very likely that that
would be the case for all $p>7$. So I computed the relation matrix
described above for products of multiweight (8,1,1,1), (7,2,1,1), (6,3,1,1),
(6,2,2,1), (5,4,1,1), (5,3,2,1), (5,2,2,2), (4,4,2,1), (4,3,3,1), (4,3,2,2) in
four even generators, and as expected they all turned out to have full rank
for $p>7$. The last multiweight to consider was (3,3,3,2). This gave the
15791$\times$8400 relation matrix mentioned above. I computed the determinants
of three randomly chosen 8400$\times$8400 submatrices, divided these
determinants out by powers of 2, 3, 5, 7, and computed the greatest common
divisor of the results. This greatest common divisor turned out to be the
product of three large primes: 529908013, 1765553401, 11899028767. I then
checked the rank of the matrix at these three primes, and it turned out to be
8399. Although there is no theoretical reason why this sort of thing should
not happen, I have never before seen anything like it. So I wrote a version of
the nilpotent quotient algorithm in \textsc{Magma}\ to handle these large
primes, and checked that four generator 5-Engel Lie algebras over these primes
have class at most 11.

\section{Locally nilpotent 5-Engel groups}

The associated Lie rings of 5-Engel groups satisfy a number of multilinear
identities in addition to the identity%
\begin{equation}
\sum_{\sigma\in\text{Sym}(5)}[a,a_{1\sigma},a_{2\sigma},a_{3\sigma}%
,a_{4\sigma},a_{5\sigma}]=0.
\end{equation}
In particular they satisfy the identity%
\begin{equation}
\sum_{\sigma\in\text{Sym}(5)}[a,a_{1\sigma},b,a_{2\sigma},a_{3\sigma
},a_{4\sigma},a_{5\sigma}]=0.
\end{equation}
To see this consider the following consequence of the group 5-Engel identity:
\[
[xy,z,z,z,z,z]=1.
\]
Since $[xy,z]=[x,z][x,z,y][y,z]$, we see that
\[
\lbrack xy,z,z,z,z,z]=[x,z,z,z,z,z][x,z,y,z,z,z,z][y,z,z,z,z,z]\text{ modulo
weight }8.
\]
It follows that in a 5-Engel group $G$,%
\[
\lbrack x,z,y,z,z,z,z]\in\gamma_{8}(G)
\]
for all $x,y,z\in G$. Setting $z=z_{1}z_{2}z_{3}z_{4}z_{5}$ and expanding we obtain%

\[%
%TCIMACRO{\dprod \limits_{\sigma\in\text{Sym}(5)}}%
%BeginExpansion
{\displaystyle\prod\limits_{\sigma\in\text{Sym}(5)}}
%EndExpansion
[x,z_{1\sigma},y,z_{2\sigma},z_{3\sigma},z_{4\sigma},z_{5\sigma}]\in\gamma
_{8}(G)
\]
for all $x,y,z_{1},z_{2},z_{3},z_{4},z_{5}\in G$, and it follows that the
associated Lie ring of $G$ satisfies the identical relation (3).

We prove Theorem 2 by showing that the associated Lie rings of 5-Engel
$p$-groups are nilpotent of class at most 10 for $p>7$. The proof follows the
same lines as our proof of Theorem 1, making use of identical relation (3)
as well as identity (2).

\end{document}